\documentclass[11pt,amsfonts, epsfig]{amsart}
\usepackage{amsmath, amscd, amssymb}
\usepackage{graphicx, psfrag} 
\usepackage{epstopdf}
\usepackage[top=1.5in, bottom=1in, left=1.2in, right=1.2in]{geometry}
\usepackage{pifont}
 \usepackage{grffile}
\usepackage{graphics}

\usepackage{float}

\usepackage{graphpap, color}
\usepackage{mathrsfs}
\usepackage{pstricks}
\usepackage{color}
\usepackage{cancel}
\usepackage[mathscr]{eucal}
\usepackage[notcite, notref]{showkeys}
\usepackage{pstricks}
\usepackage{color}
\usepackage{cancel}
\usepackage{verbatim}
\usepackage{latexsym}
\usepackage[all]{xy}

\usepackage{lipsum}
\makeatletter
\g@addto@macro{\endabstract}{\@setabstract}
\newcommand{\authorfootnotes}{\renewcommand\thefootnote{\@fnsymbol\c@footnote}}%
\makeatother

\usepackage{upgreek}

 \usepackage{wasysym}
\usepackage{ esint }

\def\bw{{\mathbf w}}
\def\virt{^{\vir}}

\numberwithin{equation}{section}

\def\sO{{\mathscr O}}
\def\sC{{\mathscr C}}

\def\sN{{\mathscr N}}
\def\sL{{\mathscr L}}

\def\sO{\mathscr{O}}

\def\sA{\mathscr{A}}

\newcommand{\CC}{\mathbb{C}}

\newcommand{\PP}{\mathbb{P}}
\newcommand{\QQ}{\mathbb{Q}}

\newcommand{\ZZ}{\mathbb{Z}}

\newcommand{\tW}{ {\widetilde{W}} }

\newcommand{\cal}{\mathcal}

\def\cL{{\cal L}}

\def\cN{{\cal N}}
\def\cO{{\cal O}}

\def\v1{{\vec{1}}}

\def\mapright#1{\,\smash{\mathop{\lra}\limits^{#1}}\,}

\def\dual{^{\vee}}

\def\sta{^\ast}

\def\virt{^{\mathrm{vir}}}

\def\sta{^{\ast}}

\def\sta{^*}

\def\lra{\longrightarrow}

\def\lsta{_{\ast}}

\newcommand{\Si}{\Sigma}

\newcommand{\si}{\sigma}

\def\begeq{\begin{equation}}
\def\endeq{\end{equation}}
\def\and{\quad{\rm and}\quad}

\def\sub{\subset}

\def\Po{{\mathbb P^1}}
\def\and{\quad\text{and}\quad}

\newtheorem{prop}{Proposition}[section]

\newtheorem{theo}[prop]{Theorem}

\newtheorem{rema}[prop]{Remark}

\newtheorem{conj}[prop]{Conjecture}

\newtheorem{defi-prop}[prop]{Definition-Proposition}
\newtheorem{defi-theo}[prop]{Definition-Theorem}

\def\sS{\mathscr S}
\def\dbar{\overline{\partial}}

\def\Ob{\cO b}

\def\sta{^\ast}

\def\sO{{\mathscr O}}

\def\sD{{\mathscr D}}

\def\beq{\begin{equation}}
\def\eeq{\end{equation}}

\def\Pf{{\PP^4}}

\def\bee{\begin{equation}}
\def\eeq{\end{equation}}

\def\sC{{\mathscr C}}

\def\bd{{\mathbf d}}

\def\barM{{\overline{M}}}

\def\mapright#1{\,\smash{\mathop{\lra}\limits^{#1}}\,}

\begin{document}

\title{On the mathematics and physics of Mixed Spin P-Fields}

\author{Huai-Liang Chang}
\address{Huai-Liang Chang,
Department of Mathematics, Hong Kong University of Science and Technology, Hong Kong} 
\email{mahlchang@ust.hk}

\author{Jun Li}
\address{Jun Li,
Shanghai Center for Mathematical Sciences, Fudan University, China; \hfil\newline 
	\indent Department of Mathematics, Stanford University,
	USA} 
\email{jli@math.stanford.edu}

\author{Wei-Ping Li}
\address{Wei-Ping Li,
Department of Mathematics, Hong Kong University of Science and Technology, Hong Kong} 
\email{mawpli@ust.hk}

\author{Chiu-Chu Melissa Liu}
\address{Chiu-Chu Melissa Liu, 
Department of Mathematics, Columbia University, USA}
\email{ccliu@math.columbia.edu}

\maketitle

\begin{abstract}  
We outline   various developments of affine and general Landau Ginzburg models in physics. 
We then describe the A-twisting and coupling to gravity in terms of Algebraic Geometry.  
We describe constructions of various path integral measures (virtual fundamental class) using the algebro-geometric technique
of cosection localization, culminating in the theory of ``Mixed Spin P (MSP) fields" developed by the authors. 
\end{abstract}


\tableofcontents

\section{Introduction}
 
In this survey we will describe several mathematical and physical theories.
\begin{enumerate}
\item The physical theory of generalized LG model (Guffin-Sharpe), and the mathematical theory of 
stable maps with P-fields and the hyperplane property in all genera (H.-L. Chang and J. Li).

\item The Fan-Jarvis-Ruan-Witten (FJRW) theory of affine LG model, and the algebro-geometric construction 
of Witten's top Chern class in the narrow case (H.-L. Chang, J. Li and W.P. Li). 

\item Witten's Gauged Linear Sigma Model (GLSM) which specializes to (1) in the Calabi-Yau (CY) phase and
(2) in the Landau-Ginzburg (LG) phase as the K\"{a}hler parameter $r\to +\infty$ and $r\to -\infty$ respectively.
The two phases are then linked via promoting the K\"{a}hler  parameter to a $\nu$-field, which leads to 
moduli spacees of Mixed Spin P (MSP) fields.
\end{enumerate}

All the above three models (stable maps with P-fields, Witten's top Chern class, MSP fields) require Kiem-Li's cosection localization. 
A ghost P-field in the CY phase can be transformed continuously to a field in the LG phase determining the spin structure on the underlying curve. 
This phonomenon is named Landau-Ginzburg transition, and is responsible for the  interaction between fields in the CY and  LG phases. 
In the end we discuss how effective the MSP moduli could be used to attack various problems, including the enumeration of
positive-genus curves in the quintic Calabi-Yau threefold. 
 
The article is written by mathematicians, aiming to include the  physics involved and the mathematics therefore stimulated, especially 
the algebro-geometric constructions by authors. 
In most sections we survey results in mathematics and in physics separately. In physics side, $C$ is a  
compact Riemann surface (worldsheet), and $K_C$ denotes its canonical line bundle.   
In mathematical  side, $C$ is an orbifold curve with markings with at worst nodal singularities, and $\omega_C$ denotes 
its dualizing sheaf. In physics part, $L$ is a $C^\infty$ complex line bundle over the compact Riemann surface $C$, 
while in mathematical  part $L$ is an algebraic line bundle over the algebraic curve $C$.
Sections of a bundle or maps between spaces, if not being mentioned ``$C^\infty$'',  are assumed to be holomorphic, i.e. 
algebraic sections/maps in algebraic geometry. Finally, we let $W_5=x_1^5+\cdots+x_5^5$ be the Fermat quintic polynomial.

\subsection*{Acknowledgments}
H.-L. Chang is partially supported by  Hong Kong GRF grant 600711 and 6301515.
J. Li is partially supported by  NSF grant DMS-1104553, DMS-1159156, and DMS-1564500.
W.-P. Li is partially supported by by Hong Kong GRF grant 602512 and 6301515.
C.-C. Liu is partially supported by  NSF grant DMS-1206667, DMS-1159416, and DMS-1564497.
This article is an expansion of J. Li's plenary talk at String Math 2015 in Sanya.

\section{Mirror Symmetry and Gromov-Witten Invariants of Quintics}

\subsection{Physics} 
A 2d supersymmetric sigma model governs maps from a fixed Riemann surface $\Sigma$ to a target manifold $X$.
When the target is a Calabi-Yau manifold, Witten \cite{MS} introduced two different ways to twist the standard supersymmetric sigma 
model, known as the A twist and the B twist, and obtained two different topological field theories, 
the A-model and B-model on $X$, denoted $A(X)$ and $B(X)$.
\begin{itemize}
\item In the A-model, the path integral over the infinite dimensional space of maps
to $X$ can be reduced to an integral over the space of holomorphic maps to $X$.
The A-model correlation functions depend on the K\"{a}hler structure but not the complex structure on $X$.
\item In  the B-model, the path integral over the infinite dimensional space of maps of $X$ can be reduced to an 
integral over the space of constant maps to $X$, i.e.,
an integral over $X$.  The B-model correlation functions depend on the complex structure but not the K\"{a}hler structure on $X$.
\end{itemize}
Given a Calabi-Yau manifold $X$, the mirror $\check{X}$ of $X$ is another Calabi-Yau manifold of the same dimension such that
\begin{equation}\label{eqn:AB}
A(X)\cong B(\check{X}),\quad B(X)\cong A(\check{X}).
\end{equation}

The expected/virtual (complex) dimension of space of holomorphic maps from a closed Riemann surface to $X$ is $\dim_{\CC}X \cdot(1-g)$, 
where $g$ is the genus of the Riemann surface. Therefore, we
expect there to be no holomorphic maps from a fixed generic Riemann surface of genus $g>1$ to $X$. By allowing the complex structure on the domain Riemann surface to vary, we obtain the 2d sigma model coupled with gravity. The A-model (resp. B-model) topological string theory on $X$ is obtained
by applying A twist (resp. B twist) to the sigma model on $X$ coupled with gravity. In the rest of this paper we will always consider
theories coupled with gravity, still denoted by $A(X)$ and $B(X)$. The mirror symmetry \eqref{eqn:AB} is still expected.
The equivalence $A(X)\cong B(\check{X})$ implies an equality of genus $g$ topological string amplitudes: 
\begin{equation}\label{eqn:FAB}
F_g^{A(X)}(q(t)) = F_g^{B(\check{X})}(t)
\end{equation}
where $t\mapsto q(t)$ is the mirror map from the moduli of complex structures on $\check{X}$ to the moduli of (complexified) K\"{a}hler classes on $X$.  The genus zero B-model
is determined by the classical variation of Hodge structures. 
In 1993, Bershadsky-Cecotti-Ooguri-Vafa (BCOV) developed the Kodaira-Spencer theory of gravity, which is a string field theory of higher genus B-model \cite{BCOV}.

A typical example is the quintic Calabi-Yau threefold $Q$, which is a degree 5 hypersurface in $\Pf$. Its mirror $\check{Q}$ is a degree 5 hypersurface in $\Pf/(\ZZ_5)^3$. 
In 1991, Candelas-de la Ossa-Green-Parkes \cite{Can} derived a formula for the 
genus zero B-model topological string amplitude $F^{B(\check{Q})}_0$ of $\check{Q}$ and the mirror map in terms
of explicit hypergeometric series, and obtained a mirror formula of the genus zero A-model topological string
amplitude $F^{A(Q)}_0$ of $Q$, which is a generating function of (virtual) numbers of rational curves in $Q$.
Mirror symmetry predictions on higher genus A-model topological string amplitudes $F_g^{A(Q)}$ (counting
genus $g$ curves in $Q$) have been obtained
by Bershadsky-Cecotti-Ooguri-Vafa at genus $g=1,2$ (\cite{BCOV}, 1993), by Katz-Klemm-Vafa at genus $g=3,4$ (\cite{KKV}, 1999), and at genus $g\leq 51$ by Huang-Klemm-Quackenbuch (\cite{HKQ}, 2007). 

Using results of BCOV \cite{BCOV} and Yamaguchi-Yau \cite{YY} and assuming mirror symmetry, 
Huang-Klemm-Quackenbush \cite{HKQ} provide an algorithm to determine $F_g^{A(Q)}(q(t))=F_g^{B(\check{Q})}(t)$  
for genus $g\leq 51$. When $g\geq 2$, the holomorphic anomaly equation determines $F_g^{A(Q)}(q)$ up to $3g-2$ unknowns. 
The degree zero Gromov-Witten invariant $N_{g,d=0}$ is known, so we are left with $3g-3$ unknowns; 
the boundary conditions at the orbifold point (which corresponds to
Landau-Ginzburg theory of the Fermat quintic polynomial in five variables) impose $\lceil\frac{3}{5}(g-1)\rceil$ contraints
on the $3g-3$ unknowns, whereas the ``gap condition'' at the conifold point imposes $2g-2$ constraints on the $3g-3$ unknowns. In summary, the holomorphic anomaly equation and the boundary conditions  determine $F_g^{A(Q)}$ 
up to $\lfloor \frac{2}{5}(g-1)\rfloor$ unknowns. When genus $g\leq 51$,  
the Gopakuma-Vafa conjecture (which relates Gromov-Witten invariants and Gopakumav-Vafa invariants)
and the Castelnuovo bound (which implies vanishing of low degree Gopakumar-Vafa invariants)
are sufficient to fix the remaining  $\lfloor \frac{2}{5}(g-1)\rfloor$ unknowns.

\subsection{Mathematics}
Gromov-Witten theory can be viewed as a mathematical theory of the A-model topological string theory. There are two approaches
to Gromov-Witten theory. Here we describe the algebro-geometric definition. For non-negative integers $d, g$, $\barM_{g}(Q,d)$ denotes  the moduli space of stable maps from genus $g$ nodal curves to $Q$ of degree $d$.  Li-Tian \cite{LT} and Behrend-Fantachi \cite{BF} construct a degree zero cycle $[\barM_{g}(Q,d)]\virt \in A_0(\barM_{g}(Q,d);\QQ)$, which is called the virtual cycle. 
Note that $\barM_g(Q,d)$ is empty when $(g,d)\in \{ (0,0), (1,0)\}$.  For $(g,d)\neq (0,0), (1,0)$, define
genus $g$, degree $d$  Gromov-Witten invariant of $Q$ by 
$$
N_{g,d}\colon =\int_{[\barM_g(Q, d)]^{vir}}1\in\mathbb Q.
$$
The genus-$g$ Gromov-Witten potential of $Q$ is given by 
$$
F^A_g(q):=
\begin{cases}
\displaystyle{\frac{5}{6}(\log q)^3+ \sum_{d=1}^\infty N_{0,d}q^d,} & g=0;\\
\displaystyle{-\frac{25}{12}\log q+ \sum_{d=1}^\infty N_{1,d}q^d,} & g=1;\\
\displaystyle{ \sum_{d=0}^\infty N_{g,d}q^d,} & g\geq 2.
\end{cases}
 $$
  
 One of the main unsolved problems in Gromov-Witten theory is to determine $F^A_g(q)$, which is a generating function of genus $g$ Gromov-Witten invariants of  $Q$.

 Using the hyperplane property in genus zero, Kontsevich \cite{Ko} proposed to use torus localization to 
calculate the genus zero Gromov-Witten invariants $N_{0,d}$.
Givental \cite{Gi} and Lian-Liu-Yau \cite{LLY} proved the mirror formula of  $F^A_0(q)$ predicted in \cite{Can}.
The BCOV mirror formula of $F^A_1(q)$ was solved in 2000's. J.Li and A. Zinger \cite{LZ} obtained a formula 
\begin{equation}\label{eqn:N-Nreduced}
N_{1,d}= N_{1,d}^{\mathrm{red}} + \frac{1}{12}N_{0,d}
\end{equation}
where $N_{1,d}^{\mathrm{red}}$ is the genus one, degree $d$ reduced GW-invariant of $Q$. 
Using \eqref{eqn:N-Nreduced} and $\mathbb C^*$-localization, Zinger proved the BCOV mirror formula  
of $F^A_1(q)$ in \cite{Zi2}. Gathmann \cite {Gath} provided an algorithm for $N_{1,d}$ using the relative GW-invariant 
formula. 

Using degeneration, Maulik and Pandharipande \cite{MP} found an algorithm which determines $N_{g,d}$ for
all genus $g$ and degree $d$: one degenerates 
the quintic of  $\Pf$ to a quartic and a $\mathbb P^3$, and than degenerates the quartic to a cubic and a $\mathbb P^3$, 
etc. In \cite[Section 0.6]{MP}, Maulik-Pandharipande described a second algorithm based on Gathmann's proposal. 
The second algorithm only requires one degeneration: one degenerates $\Pf$ to $\Pf$ and a $\mathbb P^1$ bundle over $Q$. Therefore, the second algorithm should be significantly more efficient than the first algorithm. 
Gathmann did the genus 0 and 1 cases. J. Li's degeneration formula and \cite[Theorem 1]{MP} 
(the quantum Leray-Hirsch) allow Maulik and Pandharipande to pursue Gathmann's proposal in all genera.
Maulik-Pandharipande proved that the second algorithm determines all the genus  2 invariants $N_{2,d}$ and conjectured
that it determines $N_{g,d}$ for all $g, d$; recently, L. Wu proved this conjecture in the genus 3 case.
 
We remark that the theories used by mathematicians to approach $N_{g,d}$ as above are essentially 
(1) ``hyperplane property" for $g=0,1$,  (2) torus localization formula, and (3) degeneration formula. 
They have intrinsic origin from theory of virtual cycles in mathematics. 

It remains a central problem in Gromov-Witten theory to find new effective algorithms to calculate all genus 
Gromov-Witten invariants of $Q$, with structural properties compatible with physics treatment by mirror symmetry, 
such as (quasi-)modularity of $F^A_g$ and finitely many holomorphic ambiguities with linear growth in $g$.


\section{Witten's Gauged Linear Sigma Model (GLSM)}

The same quintic polynomial $W_5=x_1^5+\ldots+x_5^5$ defines a map $\mathbb C^5\to \mathbb C$. The corresponding physical theory is the Landau-Ginzburg theory for the pair $(\CC^5,W_5)$. Since $W_5$ is invariant under the diagonal multiplicative action of $\ZZ_5$ on $\CC^5$, it descends to give an orbifold LG model $([\CC^5/\ZZ_5],W_5)$. In \cite{GLSM}, Witten embedded $Q$ into a larger background with superpotenial as follows. Let $\mathbb C^*$ act on 
$\CC^6=\CC^5\times \CC=\{(x_1, \ldots, x_5, p)\}$ with weights $(1, \ldots, 1, -5)$. 
The quotient $[\mathbb C^6/\mathbb C^*]$ has two GIT quotients:
\begin{eqnarray*}
\big((\mathbb C^5-    \{\vec 0\}  )\times \ \quad \mathbb C \ \quad \big) \ /\mathbb C^*&=&\ \ K_{\Pf}\ \ ,\\
\big(\qquad \mathbb C^5 \quad \  \, \times \,  (\mathbb C-0)\big)          \  /\mathbb C^*&=& \mathbb C^5/\mathbb Z_5.
\end{eqnarray*}
Here $K_\Pf$ is the total space of the canonical line bundle $\sO(-5)$ on $\Pf$.
The polynomial $p(x_1^5+\cdots+x_5^5)$ on $\CC^6$ is invariant under the above $\CC\sta$ action, so it 
descends to a function $\tW:[\CC^6/\CC\sta]\to \CC$.  Thus one has a picture relating generalized Landau-Ginzburg models  
\beq\label{pic1} 
 \xymatrix{& & ([\CC^6/\CC\sta],\tW) & \\
Q \  \ar @{^{(}->}[r]  &(K_\Pf,\bw)\ar @{^{(}->} [ur]^{(x_1,\cdots,x_5)\neq 0}&   & ([\CC^5/\ZZ_5],W_5) \ar @{_{(}->} [ul]_{p\neq 0}   }  
\eeq
where the restriction $\bw$ of $\tW$ to $K_{\Pf}$ is the function induced by tensoring $x_1^5+\cdots + x_5^5\in H^0(\Pf,\sO_\Pf(5))$ 
under the pairing $\sO_\Pf(-5)\otimes \sO_\Pf(5)\to \sO_\Pf$. The critical locus of the superpotential 
$\bw$ on $K_{\Pf}$ is the quintic $Q$ embedded in $K_\Pf$ 
as the subvariety defined by $p=0$ and $x_1^5+\cdots + x_5^5=0$.  The two skew arrows in Diagram \eqref{pic1} are open smooth subsets 
defined by $(x_1,\cdots,x_5)\neq 0$ and $p\neq 0$ respectively.
 
In 1993 Witten \cite{GLSM} provides a theory called  Gauged Linear Sigma Model (GLSM)
which can be considered as a sort of ``quantization" of $([\CC^6/\CC\sta],\tW)$. Here the word quantization 
means promoting variables $(x_1,\cdots,x_5)$ to fields $(\varphi_1,\cdots,\varphi_5)$ on the worldsheet (which is a connected
closed Riemann surface),  and promoting $\CC\sta$ to a principal $\CC\sta$-bundle over the worldsheet, with a gauge field. Witten's GLSM theory is parameterized by a real number $r$ called Fayet-Iliopoulos parameter, which is essentially the K\"{a}hler parameter of the symplectic quotient of $\CC^6$ by the Hamiltonian $U(1)$-action
with weights $(1,\ldots,1,-5)$.  
When  $r\to -\infty$ the GLSM is contributed by the Landau-Ginzburg model  $([\CC^5/\ZZ_5],W_5)$;
when $r\to +\infty $ the GLSM is contributed by massless instantons in $Q$, the critical locus of $(K_\Pf, \bw)$. Witten also showed that for $r\to +\infty$ the GLSM is contributed by instantons of a sigma model on the quintic threefold $Q$, i.e. holomorphic curves in $Q$. 

Witten's model suggests a few things.  Firstly, a physics theory for the Landau-Ginzburg model $(K_\Pf,\bw)$ 
should be found to undertake the specialization from GLSM to sigma model on the quintic threefold $Q$. 
Secondly, when $r>0$ is finite, contribution to GLSM is made by massive instantons in $Q$, where mass corresponds to common 
zeros of $\varphi_1,\cdots,\varphi_5$. 
In the language of mathematics, it foresees intermediate theories (obtained by varying the stability condition) other than Gromov-Witten (GW) theory of $Q$ 
(at $r\to +\infty$) or Landau-Ginzburg (LG) theory of $([\CC^5/\ZZ_5],W_5)$  (at $r\to -\infty$).


When \cite{GLSM} first appeared, it was far from a mathematical theory:
as Witten's \cite{GLSM} is a gauged field theory defined by using path integral, mathematicians need to find conditions to ensure
the convergence of the path integral and then substitute the infinite dimensional path integral measure with certain finite dimensional construction to
obtain rigorous mathematical definitions. Moreover, theories for every spaces in \eqref{pic1} were neither twisted nor coupled to gravity. 
  

  Witten's GLSM tells us, once mathematicians can possibly achieve finite dimensional constructions which lead to rigorous mathematical definitions,   
the (massive or massless) theories for  $Q$, $(K_\Pf,\bw)$,  $([\CC^5/\ZZ_5],W_5)$, and possibly even the universal $([\CC^6/\CC\sta],\tW)$, should determine each other; namely, they should be 
``equivalent" theories. However, what are the explicit relations among amplitudes from any two different theories?   
Could these conjectural equivalences help determine all of them, or just reduce three sorts of mysteries to one that is still mysterious? 
In later sections of  this paper, we will discuss solutions to the above questions on constructions and relations. 

In the following we shall call $([\CC^5/\ZZ_5],W)$ an affine Landau-Ginzburg model as $\CC^5$ is an affine space,  to be distiguished 
from $(K_{\PP^4},\bw)$, which is a general Landau-Ginzburg model. 

\section{Hyperplane Property, Ghost, and P-field} \label{sec:hyperplane}

\subsection{Physics: Guffin and Sharpe}  Guffin and Sharpe consider the A-twisting of the 
LG model $(K_\Pf,\bw)$ suggested by Witten, and show its amplitudes  are equal to genus zero GW invariants
of $Q$. This can be viewed as the hyperplane property in physics, at least in genus zero. 
To describe the matter field $\varphi:C ``\longrightarrow" K_\Pf$ (where $\varphi$ and $C$ are smooth), 
in additional to classical $\varphi_1,\cdots,\varphi_5$ as $C^\infty$ sections of a line bundle $L$ 
on $C$, the noncompact direction needs to be twisted by the canonical line bundle $K_C$ of $C$ and 
considered as 
\beq\label{GSP}
p\in C^\infty (C,K_C\otimes L^{\otimes -5}),
\eeq
so  a term in the Lagrangian becomes a top form on $C$ and can be integrated to make sense of the action.  
Guffin and Sharpe showed that , in the genus zero case, their integral (as an invariant associated to enumerating 
curves mapped to $(K_\Pf,\bw)$) is equal to
$$
e(E_d)\cap [X_d] \in H_0(X;\QQ)=\QQ
$$
where $e(E_d)$ is the Euler class of the finite rank complex vector bundle
$$ 
E_d=\underset{\varphi\in X_d}{\bigcup} H^0(C,\varphi\sta\sO_\Pf(5))
$$
over the smooth compact complex orbifold $X_d=\{\varphi:C\to \Pf \mid \deg\varphi=d, g(C)=0\}$.
By Kontsevich's hyperplane property,  $e(E_d)\cap [X_d]$ is the genus zero, degree $d$ GW invariant $N_{g=0,d}$  of $Q$.

\subsection{Mathematics: Hyperplane Property}  We explain the hyperplane problem in mathematics.
Fixing the degree $d$, for each genus $g$, over the finite dimensional compact space
$$
Y_{g,d}:=\{\varphi:C\to \Pf\ \textup{ stable map }\mid \deg\varphi=d, g(C)=g\}
$$ 
there are two unions of vector spaces 
$$ 
V_{g,d}:=\bigcup_{\varphi\in Y_g} H^0(C,\varphi\sta\sO_\Pf(5))
    \and
V'_{g,d}:=\bigcup_{\varphi\in Y_g} H^1(C,\varphi\sta\sO_\Pf(5)),
$$
defining two sheaves over $Y_{g,d}$. 
By Riemann-Roch formula, the difference of the dimensions $\dim V_{g,d}|_{\varphi}-\dim V'_{g,d}|_{\varphi}$ is $5d+1-g$, which
is independent of $\varphi\in Y_{g,d}$.

If $g=0$, one can show $V'_{0,d}=0$.  Thus $V_{0,d}$ has constant dimensional fiber over $Y_{0,d}$, namely 
$V_{0,d}$ is a complex vector bundle of rank $5d+1$ over $Y_{0,d}$.  The hyperplane property of Kontseviech says 
\beq\label{hyper}
N_{0,d}= \deg\big(e(V_{0,d})\cap [Y_{0,d}] \big)\in \QQ.
\eeq 
where $\deg: A_0(X_d;\QQ)\to \QQ$. This reconstructs the enumeration of rational curves in $Q$  from the information of $\Pf$. 
The identity \eqref{hyper} is an easy consequence of virtual cycle theory \cite{KKP}. By \eqref{hyper}, one may compute $N_{0,d}$
by torus localization as $V_{0,d}$ and $Y_{0,d}$ admit a $(\CC^*)^4$-action inherited from that on $\Pf$. 
  
When $g>0$, everything above fails unfortunately. For example when $g=1$, $Y_{g=1,d}$ contains essentially two kinds of components. 
The two components collect maps of different forms (see Figure 1 below).

 
\begin{figure}[h] 
\begin{center}
 \includegraphics[scale=0.4]{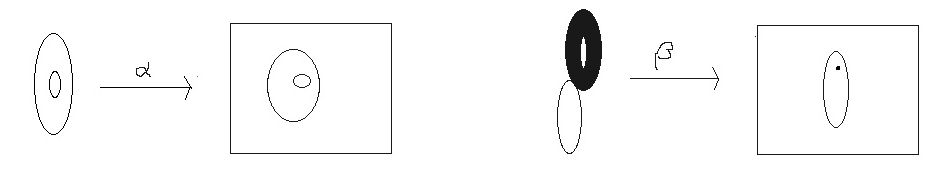}
\end{center}
\caption{graphs for an honest map $\alpha$ and a ghost map $\beta$.}
\end{figure}
The main component of $Y_{1,d}$ consists of maps looking like $\alpha$, which has positive degree on the genus one component of the curve. 
The other component of $Y_{1,d}$ consists of maps looking like $\beta$, which contracts the (black) genus one component to one point, 
and has positive degree on the (blank) genus zero  $\Po$. 
For the curve and maps being indicated as $\alpha$, 
every element in $H^1(C,L^{\otimes 5})=H^0(C,K_C\otimes L^{\otimes -5})\dual$  vanishes since $L=\alpha\sta\sO(1)$ has positive degree and  $K_C\otimes L^{\otimes -5}$ has negative degree. Thus the P-field must vanish for $\alpha$ (recall that $P\in H^0(C,K_C\otimes L^{\otimes -5})$),
or equivalently, $V'_{1,d}|_\alpha=0$. However, for $\beta$ one can have nonzero one form on the elliptic component, which extends
by zero to give a section of $K_C\otimes L^{\otimes -5}$. This corresponds to the fact that $\beta$ contracts a genus one component of the
curve $C$, for which reason we say $\beta$ is a ghost map. The black genus one component is a
``ghost'' component, and \textit{P-field can survive on a ghost}.  One easily calculates $V_{1,d}'|_\beta\cong \CC$.  
  
   We remark that in the approach determing $N_{1,d}$ in \cite{LZ,Zi}, a key issue is to locate the contribution of the ghost in the counting. 
In their formula $N_{1,d}=N_{1,d}^{\mathrm{red}}+\frac{1}{12}N_{0,d}$,  the term 
$N_{1,d}^{\mathrm{red}}$ is the contribution from maps of type $\alpha$, and $\frac{1}{12}N_{0,d}$ 
is the contribution from maps of type $\beta$, where the $\frac{1}{12}$ 
comes from integrating out all the P-fields living on the ghost (black) elliptic component of $\beta$. 
For our ultimate purpose to approach $N_{g,d}$ for larger $g$, locating the contribution of P-field (including ghosts) 
becomes very difficult and out of control.

Since $V'_{1,d}|_\alpha=0$ and $V'_{1,d}|_\beta\cong \CC$,  $V_{1,d}'$ has fiber rank jumping over 
$X_{1,d}$, and by Riemann-Roch $V_{1,d}$ also does and is not 
a vector bundle over $Y_{1,d}$.  As  the Euler class is only defined for vector bundles, 
$e(V_{1,d})$ no longer makes sense.  
It is natural to ask how the hyperplane property \eqref{hyper} should be modified, so that the information of $\Pf$ 
can be used to reconstruct enumeration of higher genus curves in $Q$ in mathematics (namely only finite dimensional construction allowed).
  
 After A-twisting, the topological string theory (with supersymmetry) admits a mathematical counterpart 
called ``virtual cycle"\footnote{also called virtual fundamendal class}.
 As virtual cycle (\cite{LT}) is governed by tangent-obstruction (deformation) theory (in physics words, after A-twisting, 
the zero mode of fermions over SUSY fixed loci, even if the loci is singular, recovers the path integral algebraically), we
may view the above problem of higher genus hyperplane property in the following way. Let $f:C\to Q$ be a point in
$$
f=[f_1,\cdots,f_5]\in \barM_g(Q,d):=X_{g,d} \sub \barM_g(\Pf,d)=Y_{g,d}.
$$ 
The exact sequence $0\to T_Q\to T_\Pf|_Q\to \sO(5)|_Q\to 0$ induces the following long exact sequence
\beq\label{seq0}\begin{CD}
0@>>> H^0(C,f\sta T_Q) @>>> H^0(C,f\sta T_\Pf) @>>> H^0(C,f\sta \sO(5))\\
@>>>H^1(C,f\sta T_Q) @>>> H^1(C,f\sta T_\Pf) @>>> H^1(C,f\sta \sO(5))@>>>0.
\end{CD}\eeq
Every vector space above has geometric meanings, namely the sequence \eqref{seq0} is identical to
\beq\label{seq1}\begin{CD}
0@>>>  T_{f,X_{g,d}}@>>> T_{f,Y_{g,d}}@>>>\Ob_{f, X_{g,d}/Y_{g,d}}\\
@>>>\Ob_{f,X_{g,d}} @>>> \Ob_{f,Y_{g,d}} @>>> \Ob^{\mathrm{higher}}_{f,X_{g,d}/Y_{g,d}}@>>>0.
\end{CD}\eeq
where
\begin{itemize}
\item $T_{f,X_{g,d}}$ and $T_{f,Y_{g,d}}$ are the first order deformations of $f$ in $X_{g,d}$ and $Y_{g,d}$ 
(relative to moduli space of genus $g$ nodal curves) respectively;
\item $\Ob_{f,X_{g,d}}$ and $\Ob_{f,Y_{g,d}}$ are obstructions to deforming $f$ in $X_{g,d}$ and $Y_{g,d}$ respectively;
\item $H^0(C,f\sta \sO(5))$, which contains the element $f_1^5+\cdots+f_5^5$, 
is the obstruction\footnote{ because $f_1^5+\cdots+f_5^5=0$ characterises the $f\in Y_{g,d}$ lies in $X_{g,d}$} 
of an element in $Y_{g,d}$ to be in $X_{g,d}$, namely the relative obstruction $\Ob_{f,X_{g,d}/Y_{g,d}}$;
\item $H^1(C,f\sta \sO(5))$ is the \textit{higher obstruction} of a point in $Y_{g,d}$ to lie in  $X_{g,d}$. 
\end{itemize}

 Now recall that the tangent and obstruction theory would determine the virtual cycle  (path integral measure),  and the two terms in the left column  in \eqref{seq1}  are tangent and obstruction theories of $X_{g,d}$, therefore are responsible for the Gromov-Witten invariant 
$N_{g,d}$ of the quintic Calabi-Yau threefold $Q$. The two terms in the middle column are tangent and obstruction theories 
of $Y_{g,d}$ which parametrizes maps to $\Pf$. To solve the hyperplane property problem, 
one should combine the right column in \eqref{seq1} with the datum of $Y_{g,d}$.
If this can be done then one may expect to recover $N_{g,d}$.
 
We observe that the last term (higher obstruction) $H^1(C,f\sta\sO(5))$ is dual to the space $H^0(C,K_C\otimes f\sta\sO(-5))$ 
of algebraic P-fields (c.f. \eqref{GSP}), which we may add it to the moduli space $Y_{g,d}=\barM_g(\Pf,d)$
of stable maps to $\Pf$ to form the moduli of stable maps to $\Pf$ with 
P-fields\footnote{now allow $f:C\to Q$ to be more general $\varphi:C\to\Pf$}
\beq\label{P}
Y_{g,d}^p:=\barM_g(\Pf,d)^p=\{[\varphi:C\to \Pf]\in \barM_g(\Pf,d)),  \rho \in \Gamma(C,K_C\otimes \varphi\sta \sO(-5))\},
\eeq
where  $\rho$ is called an algebraic ``P-field" as its analogue in \eqref{GSP}.
As the obstructions to deforming $\varphi$ and $\rho$ lie in $H^1(C,\phi\sta T_\Pf)$ and  
$H^1(C,K_C\otimes \varphi\sta(\sO(-5)))=H^0(C,\varphi\sta\sO(5))\dual$ respectively, the deformation theory of $Y_{g,d}^p$ 
is given by the middle and right columns in \eqref{seq0}. If one is able to define a virtual cycle for $\barM_g(\Pf,d)^p$, then 
it is expected to be ``equivalent" to the virtual cycle of $X_{g,d}$, and the hyperplane problem is solved.
Now the difficulty appears because $\barM_g(\Pf,d)^p$ is {\em non-compact} due to the presence of $P$-fields: for example, over the ghost map $\beta$, 
the P-field can be any element in $\CC$ that is unbounded. This difficulty is then overcome by the invention of ``cosection localization" 
by Y.H. Kiem and Jun Li, along with H.L. Chang's observation that 
``{\it the supersymmetry variation of the superpotential on worldsheet defines a cosection, which solves Witten's equation in the general Landau-Ginzburg theory}."

 We now describe the algebro-goemetric results of H.-L. Chang and J. Li \cite{CL} discovered based on the above reasoning. 
In the definition of \eqref{P},  the data $([\varphi,C],\rho)$ is equivalent to the data $ (C,L,\varphi_1,\cdots,\varphi_5, \rho)$, since  the map 
$\varphi$ is equivalent to the line bundle $L=\varphi^*\sO_{\Pf}(1)$ with five sections $(\varphi_1,\cdots,\varphi_5)$ of $L$. 
We regard   $\barM_g(\Pf,d)^p$  as a space  of ``maps from curve to $K_\Pf$." 
The moduli stack $\barM_g(\Pf,d)^p$ has a perfect obstruction theory
relative to the smooth Artin stack $\sD=\{(C, L)\}$. 
At $\xi=[(C, L, \varphi_i,\rho)]\in \barM_g(\Pf,d)^p$, the (relative) obstruction space of deforming $\xi$ is 
$$
\sO b_{M/\sD}|_\xi=H^1(L)^{\oplus 5}\oplus H^1(L^{\vee 5}\otimes \omega_C).
$$
There exists a cosection 
$$
\sigma\colon \sO b_{M/\sD} \to \sO_{\barM_g(\Pf,d)^p}
$$
constructed as follows. Let 
$$
(\dot\varphi_1,\ldots, \dot\varphi_5, \dot\rho)\in H^1(L)^{\oplus 5}\oplus H^1(L^{\vee 5}\otimes \omega_C)=\sO b_{M/\sD}|_\xi.
$$
Define
$$\sigma|_{\xi}(\dot\varphi_1,\ldots, \dot\varphi_5,\dot\rho)\colon =\dot\rho\sum_{i=1}^5\varphi_i^5+\rho\sum_{i=1}^55\varphi_i^4\dot\varphi_i.$$

The degeneracy locus $D(\sigma)$ of the cosection $\si$ consists of $\xi$ such that $\sigma|_{\xi}$ is zero, i.e., $\sigma|_{\xi}(\dot\varphi_1,\ldots, \dot\varphi_5,\dot\rho)=0$ for all $\dot\varphi_i$ and $\dot\rho$.   Thus 
$$D(\sigma)=\{\xi \in \barM_g(\Pf,d)^p\, |\, \rho=0\hbox{ and } \sum_{i=1}^5\varphi_i^5=0\}=\barM_g(Q,d)\subset \barM_g(\Pf,d).$$
This corresponds to the fixed loci of supersymmetry (SUSY) in path integral.
The expression of the cosection $\sigma$ comes from supersymmetry variation $\delta$ (in physics) 
applied to $p\cdot W_5=p(x_1^5+\ldots+x_5^5)$  where $p$ and $x_i$ live on the worldsheet, via H.L. Chang's observation.
 
 Since  $\rho$ is a section, the moduli space $ \barM_g(\Pf,d)^p$ is not proper (when $g\ge 1$) and hence cannot be used to define invariants. However, the degeneracy locus $D(\sigma)$ is the moduli space 
$\barM_g(Q,d)$ of stable maps to the quintic threefold $Q$ and thus proper. Using cosection localization developed by  Y.H. Kiem and J. Li \cite{KL}, H.L. Chang and J. Li
 constructed \cite{CL} the cosection localized virtual cycle for Landau-Ginzburg theory
$$
[ \barM_g(\Pf,d)^p]^{\mathrm{vir}}_{\mathrm{loc}}\in A_*D(\sigma)=A_*\barM_g(Q,d).
$$
 As always one defines  the $P$-fields GW invariants 
$$
N^p_{g, d}=\int_{[\barM_g(\Pf,d)^p]^{vir}_{loc}}1\in \mathbb Q.
$$
H.L Chang and J. Li proved the following.
 \begin{theo}[H.L. Chang - J. Li \cite{CL}]\label{GW-GSW} 
The GW invariant of the quintic threefold $Q$ equals the P-fields GW-invariant up to a sign:
$$
N_{g,d}=(-1)^{d+g+1}N^p_{g, d}.
$$
\end{theo}
 
The advantage of this result is that $F^A_g(q)=\sum_d N_{g, d}q^d$ now becomes the amplitude of a theory valued in $K_{\Pf}=\big((\mathbb C^5-\vec 0)\times \mathbb C\big)/\mathbb C^*$.

In conclusion, the invariant enumerating maps from curves to $(K_\Pf,\bw)$ is equal to 
the invariant enumerating maps from curves to $Q$, up to a sign. This generalizes the genus zero case to the positive genus case, 
and solves the hyperplane property problem.

\section{Fields Valued in Two GIT Quotients}
 
\subsection{Physics: GLSM}

 From \S \ref{sec:hyperplane}, we see the curve-enumerating theories  for $Q$ and $(K_\Pf,\bw)$ in Diagram \eqref{pic1} are both established in mathematics for all genera. 
Witten's \cite{GLSM} suggests the theory for $([\CC^5/\ZZ_5,W_5])$ at the lower right corner of  \eqref{pic1} should also exist, and match the physical theory of A-twisted LG model coupled with gravity.

\subsection{Mathematics}
 We now consider the space of maps from curves to each target in \eqref{pic1}, viewed as a sort of ``quantization" of \eqref{pic1}. 

 The previous sections tells us the space $\barM_g(\Pf,d)^p$ of all  ``maps to $(K_\Pf,\bw)$" is the set of all  $(C,L,\varphi,\rho)$
 where   $\varphi=(\varphi_1,\cdots,\varphi_5)$ is section of $L^{\oplus 5}$, $\rho$ is a section of $\omega_C\otimes L^{\otimes -5}$, and 
\[
 \varphi= (\varphi_1,\cdots,\varphi_5) \ \text{ has  no  zeros on }\ C \label{1} \tag{+}
\]
 so that $\varphi_i's$ define an honest map to $\Pf$. 
 Without  the condition \eqref{1}, one obtains  a huge Artin stack $\sA rt$ of all $(C,L,\varphi,\rho)$ for arbitrary 
$\varphi\in\Gamma(C,L^{\oplus 5})$ and $\rho\in\Gamma(C,\omega_C\otimes L^{\otimes -5})$. 
The stack $\sA rt=\{(C,L,\varphi,\rho)\}$  should be viewed as the moduli space of maps to $[\CC^6/\CC\sta]_{(1,1,1,1,1,-5)}$.
$\barM_g(\Pf,d)^p$  is the open substack of objects in $\sA rt$ subject to condition \eqref{1}, which
corresponds to the open substack $K_\Pf\sub [\CC^6/\CC\sta]$ in \eqref{pic1}  
defined by  $(x_1,\cdots,x_5)\neq (0,\cdots,0)$. After quantizing it translates to the requirement \eqref{1}, as $\varphi_1,\cdots,\varphi_5$ are the five fields  
promoted from the five coordinates $x_1,\cdots,x_5$.

Parallelly, since the open substack $[\CC^5/\ZZ_5]\sub [\CC^6/\CC\sta]$  is defined by $p\neq 0$ in \eqref{pic1}, to define a theory whose target is 
$([\CC^5/\ZZ_5],W_5)$, one analogously expects to pick up the open substack of $\sA rt$ subject to the condition 
 \[
\rho \ \text{ has  no  zeros on }\ C   \label{2} \tag{--}
\]
as $\rho$ is the field promoted from coordinate $p$ in \eqref{pic1}. 
Namely $\rho$ trivializes $\omega_C\otimes L^{-5}$, or equivalently, gives an isomorphism $L^{\otimes 5} \mapright{\cong} \omega_C$. 
One then expects the theory  of $([\CC^5/\ZZ_5],W_5)$ to start with the moduli space of all $(C,L,\varphi)$ where
\begin{enumerate}
\item $L$ is a fifth root\footnote{sometimes called $5$-spin structure on $C$;} of $\omega_C$,  and
\item $\varphi=(\varphi_1,\cdots,\varphi_5)$ is an arbitrary section of $L^{\oplus 5}$. 
\end{enumerate} 
We denote this moduli as $\barM_g^{1/5,5p}$ where $1/5$ denotes the 5-spin structure and $5p$ indicates that an object consists of
five sections $\varphi_1,\cdots,\varphi_5$ of $L$ (by abuse of notation).
 
 When one quantizes  every space in \eqref{pic1}, one then obtains two open substacks (subspaces) of the common huge Artin stack as follows:
  \beq\label{pic2} 
 \xymatrix{& & \{(C,L,\varphi,\rho)\} & \\
\barM_g(Q,d) \  \ar @{^{(}->}[r]  &\barM_{g}(\Pf,d)^p  \ar @{^{(}->} [ur]^{\varphi\, \text{nowhere} \, 0}&   &  \barM_g^{1/5,5p}   \ar @{_{(}->} [ul]_{\rho\,\text{nowhere}\, 0}   }  
\eeq

Naturally one wonders whether the substack $\barM_g^{1/5,5p}$ at the bottom right corner has a virtual cycle, with which 
intersections represent invariants of the Landau-Ginzburg model $([\CC^5/\ZZ_5],W_5)$ 
from physics, as Witten predicted. 
Coincidentally, around 2010 H. Fan, T. Jarvis, and Y. Ruan carried out a construction of an A-side theory of which 
the target may be any affine LG space $([\CC^n/G],\bw)$, where $G$ is a finite group and 
the ``superpotential'' $\bw$ is a $G$-invariant polynomial on $\CC^n$. Their approach to the affine Landau-Ginzburg model $([\CC^n/G],\bw)$ originates from a different line 
in history, namely the gauged WZW model, Witten equation, and Hamiltonian Floer theory,  which we brief in \S \ref{sec:affineLG} below.

 
\section{Affine LG Phase and Spin Structure}  \label{sec:affineLG}
 
\subsection{Physics: SUSY A-twisted LG Theory Coupled To Gravity} 
The classical Landau-Ginzburg theory on the A-side follows a different line of development in history.
In \cite{Mg} E. Witten  conjectured that descendant integrals on moduli spaces
of stable curves $\barM_{g,n}$ satisfy the KdV equations, and the string equation (proved by Witten)
and the KdV equations uniquely determine all descendant integrals from the initial
value $\int_{\barM_{0,3}}1=1$.  Witten's conjecture was first proved by 
Kontsevich \cite{Kont1} by stratification of $\barM_{g,n}$ and matrix model. 
For the purpose to generalize above to $N$ matrix model, Witten in \cite{N} considered
A-twisted  gauged WZW model targeting $SU(2)/U(1)$ coupled with gravity, and obtained a topological theory 
which he conjectured \cite{GWZW} (refining/twisting the minimal model of \cite{KLi} et.al.) to solve 
the generalized KdV hierarchies ($N$-matrix model). 
 
Witten's A-twisted theory is localized to the SUSY fixed locus consisting of objects 
almost definable in algebraic geometry.
Let $\barM_g^{1/r}$ denote the moduli space of Riemann surfaces $C$ (with at worst nodal singularities)  together with a line bundle $L$ such that $L^{\otimes r}\cong K_C$.
Mathematically $\xi=(C,L)$ is referred as a $r$-spin curve. One may also add orbifold marked points on $C$ but we omit them here for simplicity. 
Witten roughly argued that $\barM_g^{1/r}$ is smooth and compact. He set 
$$
M^\infty=\bigcup_{(C,L)\in \barM_g^{1/r}}C^\infty(C,L)
$$ 
where $C^\infty(C,L)$ is the space of $C^\infty$ sections $u$ of $L$. 
  
The topological correlation function of Witten's theory  amounts to counting the intersection number of the zero section of the (infinite rank) bundle
\beq\label{spinMQ}
E^{\infty}=\underset{(\xi,u)\in M^\infty }\bigcup\Omega^{(0,1)}_C(L)\lra  M^\infty 
\eeq
with the graph of the section
\beq\label{wi}
(\xi=(C,L),u)\mapsto   s_W(\xi,u):=\dbar u+ r (\bar u)^{r-1}
\eeq
and possibly with insertions (\cite{GWZW}) such as  gravitational descendents 
(if one adds markings on each $C$).
Note that we may choose a K\"{a}hler metric on the Riemann surface $C$ and a Hermitian metric on the line bundle $L$,
so that  $(\bar u)^{r-1}$ becomes a section of $(\bar L)^{\otimes (r-1)}\cong L^{\otimes (1-r)}\cong \overline{K}_C\otimes L$, where $\dbar u$ lives.  In short the theory counts solutions of 
\beq\label{count}
\dbar u+ r (\bar u)^{r-1}=0.
\eeq
The Euler class of $E^\infty$ localized by Witten's section $s_W$ is then called ``Witten's top Chern class", a core object in the definition of the theory.
 
For the purpose to interprete Witten's correlation function more directly,  one may regard it 
as the A-twisted (and coupled to gravity) version 
of the   ``Landau-Ginzburg theories"  
defined in \cite{Vafa}, \cite{Ito} (also c.f. \cite{Ce}). Vafa, et.al.'s model build the 
Landau-Ginzburg structure directly in the Lagrangian. Namely, it is a path integral whose configuration space of fields is the set of maps 
from the worldsheet to $([\CC^n/G],\bw)$, with fermions coupled with terms as 
$$
s_\bw=( \dbar u_i+\overline{\partial_{u_i}\bw(u_1,\cdots,u_n)})_{i=1}^n,
$$
 and the contribution to the theory comes from solutions of 
\beq\label{Weqn}
\dbar u_i+\overline{\partial_{u_i}\bw(u_1,\cdots,u_n)}=0\qquad \text{for all} \ i=1,\cdots,n
\eeq
generalizing \eqref{count} where $n=1,\bw=x^r$.

However, the theories in  \cite{Ito}, \cite{Vafa}  are not coupled with gravity,
and the group $G$ is trivial $G=\{e\}$. 
It was then later understood  (by Fan-Jarvis-Ruan etc.)  that Witten's model is using 
$([\CC/\ZZ_r],\bw=x^r)$, whose state spaces are indexed by the monodromy weights of the $r$-spin bundle at   markings.  \black

 \subsection{Mathematics: FJRW Invariants}  

 Based on Witten's infinite dimensional Euler class model (with section to be $s_\bw$), 
Fan-Jarvis-Ruan  \cite{FJR1, FJR2} used analytic methods to construct the Witten's top Chern class,
and defined correlators of a Cohomological Field Thery (CohFT) by capping the Witten's top Chern class
with states of the Landau-Ginzburg model $([\CC^n/G],\bw)$. 
Fan-Jarvis-Ruan's pioneer work is now known as FJRW invariants associated to the singularity  $([\CC^n/G],\bw)$. 
FJRW invariants of special $ADE$ type singularities can be enumerated and are governed by 
the Kac-Wakimoto/Drinfeld-Sokolov hierarchies \cite{LRZ}, generalizing \cite{FSZ}'s proof of Witten's $r$-spin conjecture.

In FJRW theory, Witten's top Chern class is constructed in differential geometry via perturbing \eqref{Weqn}.
It can also be constructed in algebraic geometry  without pertburbing \eqref{Weqn}. 
The algebro-geometric constructions (in the narrow case) were carried out by Polishchuk-Vaintrob \cite{PV}, by Chiodo  \cite{Chi}, 
and by H.L. Chang, J. Li and W.P. Li \cite{CLL}. For our purpose to provide a field theory valued in $([\mathbb C^5/\mathbb Z_5],W_5)$, we brief the construction in \cite{CLL} here, using the version with  markings. 
Recall that the moduli  $\barM_g^{1/5,5p}$ for $([\CC^5/\ZZ_5],W_5)$ requires a fifth root $L$ of $\omega_C$, which 
does not exist if $\deg \omega_C=2g-2$ is not divisible by five. 
One thus extends the setup by allowing $C$ to be a twisted curve with markings (which can be scheme points or stacky points). 
Thus our field valued in $\big( \mathbb C^5  \times (\mathbb C-0)\big)/\mathbb C^*=[\mathbb C^5/\mathbb Z_5]$ consists of 
$$
\xi=(\Sigma^C,  C, L, \varphi_1,\ldots,\varphi_5, \rho)
$$
where $(\Sigma^C, C)$ is a pointed twisted curve with markings $\Sigma^C$ possibly stacky, $L$ is an invertible sheaf on $C$, 
$\varphi_i\in H^0(L)$, and $\rho\in H^0(L^{\vee 5}\otimes \omega^{\log}_C)$ with $\omega^{\log}_C=\omega_C(\Sigma^C)$, 
and the corresponding property of \eqref{2}
 \[
\text{ the section}\ \rho \ \text{is  nowhere vanishing}   \label{3} \tag{*}
\]
is required. This implies $L^{\vee 5}\otimes \omega^{\log}_C\cong \sO_C$, or equivalently $L^{\otimes 5}\cong \omega^{\log}_{C}$. Therefore $(\Sigma^{C},C, L)$ is a $5$-spin curve and $(\varphi_1, \ldots, \varphi_5)$ gives five fields. 
We get a moduli space of $5$-spin curves with five fields:
\begin{eqnarray*}
\barM_{g,\gamma}^{1/5,5p}=\{(\Sigma^{C}, C, \cL, \varphi_1,\cdots,\varphi_5, \rho) \mid \text{$\rho$ is nowhere zero}\}.
\end{eqnarray*}
Here $\gamma$ is the monodromy data: if $\Sigma_j$ is a stacky marking on $C$, then $\mu_5$ acts on $L|_{\Sigma_j}$ with weight 
$\gamma_j={ \exp}(2\pi i r_j/5)$ where $1\le r_j\le 4$ and we call $\gamma_j$ narrow. If $\Sigma_j$ is a scheme marking, we call it broad, and it corresponds to $\gamma_j=1$.  

Similar to the case $\barM_{g}(\Pf,d)^p$, the moduli stack $\barM_{g,\gamma}^{1/5,5p}$
has a perfect obstruction theory relative to the smooth Artin stack $\sD=\{(\Sigma^{C}, C, L)\}$. 
There exists a cosection $\sigma\colon \sO b\to \sO_{\barM_{g,\gamma}^{1/5,5p}}$ whose degeneracy locus  is 
$$
D(\sigma)=\{\xi\in \barM_{g,\gamma}^{1/5,5p} \,|\, \varphi_i=0 \hbox{ for all $i$}\}=\barM_{g,\gamma}^{1/5}=\{(\Sigma^{C}, C, L)\, |\, L^{\otimes 5}\cong \omega^{\log}_C\},
$$
which is the moduli space of $5$-spin curves. 
\begin{theo}[H.L. Chang - J. Li - W. P. Li \cite{CLL}] 
The (narrow) FJRW invariants can be constructed using cosection localized virtual cycles of $\barM_{g,\gamma}^{1/5,5p}$:
$$
[\barM_{g,\gamma}^{1/5,5p}]^{\mathrm{vir}}_{\mathrm{loc}}\in A_*\barM_{g,\gamma}^{1/5,5p}.
$$
\end{theo}

The Witten equation mentioned in \eqref{Weqn}, in this case, becomes
\begin{eqnarray}\label{Witten's-eq}
\bar\partial s_i+\overline{\partial_{x_i}W_5(s_1, \ldots, s_5)}=0, \quad i.e., \quad \bar\partial s_i+5\overline {s_i^4}=0 \qquad \text{for}\ i=1,\cdots,5.
\end{eqnarray}
This is used to construct Witten's top Chern class to define invariants on the moduli space of $5$-spin curves. 
From Witten's equation \eqref{Witten's-eq}, 
the term $\bar \partial s_i$ gives the obstruction to extending a holomorphic section. 
Thus the left hand side of \eqref{Witten's-eq} is a $C^\infty$ section of the obstruction sheaf of the moduli of spin 
curves with fields. Substituting the complex conjugate in the Witten's equation by the Serre duality, the left hand side of 
\eqref{Witten's-eq}  becomes a smooth inverse of cosection. 
The Mathai-Quillen setup in \eqref{spinMQ} and   \eqref{wi} generalize naturally here and the Euler class localized near solution of Witten equations would be equal to the Kiem-Li's virtual cycle localized via cosection $\sigma$. 

 We remark here that the form \eqref{Witten's-eq} indicates the virtual cycle of $([\CC^5/\ZZ_5],W_5)$ is the five self-intersection of the virtual cycle of $([\CC/\ZZ_5],x^5)$ (each defined by using $\dbar s+5\overline{s^5}=0$ which is 
 \eqref{count} for the case $r=5$). This remarkable property  is related to self-tensor product of conformal field theories and 
is discussed in \cite{FJR1} or \cite{CR}.


There is an important subclass of FJRW invariants: those with insertions $-\displaystyle\frac{2}{5}$. Let $C$ have $k$ markings with all $\gamma_j=\zeta^2$ for $1\le j\le k$ where 
$\zeta={\exp}(2\pi i /5)$. For simplicity we write $\gamma=(\gamma_j)_{j=1}^k=2^k$. Define
$$
\Theta_{g, k}\colon =\int_{[\barM_{g,2^k}^{1/5,5p}]^\mathrm{vir}_\mathrm{loc}}1\in\mathbb Q, \quad \hbox{for $k+2-2g=0$ {mod} $5$}.
$$
 It is shown \cite{CLLL2} that $\{\Theta_{g,k}\}_{g,k}$ determine all FJRW invariants with descendents for  the quintic LG space  $([\CC^5/\ZZ_5],W_5)$, where an explicit formula will be given in \cite{twFJRW}. 
For this reason we call $\{\Theta_{g,k}\}_{g,k}$ the primary FJRW invariants.

\section{The Puzzle to Link Invariants in Opposite Phases}
\subsection{Mathematics}

We have seen that the three moduli spaces at the bottom of Diagram \eqref{pic2} admit virtual fundamental classes, 
while the moduli space $\sA rt:= \{(C,L,\varphi,\rho)\}$ at the top of \eqref{pic2} does not,  because $\sA rt$ is not a Deligne Mumford stack.  
One can introduce all possible stability conditions to define open substacks of $\sA rt$ that are Deligne Mumford, just as $\barM_{g}(\Pf,d)^p$   
and $\barM_g^{1/5,5p} $, and then construct virtual classes (path integral measures) for them as we defined 
$[\barM_{g}(\Pf,d)^p]_{\mathrm{loc}}\virt$ and $[\barM_g^{1/5,5p}]\virt_{\mathrm{loc}}$ using P-fields and cosection machinery. 
This is the step that most groups are taking. The theory of $\epsilon$-stability and quasimaps (\cite{FK1} \cite{StableQ}) are developed, for example.

 On the other hand, introducing new stability conditions means there are invariants other than the original  $N_{g,d}$'s.. 
Whether these new invariants (defined by new stability conditions) can simplify enumeration of $N_{g,d}$'s
 or give structures for $N_{g,d}$ predicted by the B-side, is not easy at all. 
Following Witten's GLSM, we wishfully expect knowing FJRW invariants $\Theta_{g,k}$'s of $(\CC^5/\ZZ_5],W_5)$ 
would help us to understand/enumerate GW invariants $N_{g,d}$'s of $Q$. We would like to 
know whether, and how exactly, the invariants $N_{g,d}^p$'s defined by $[\barM_g(\Pf,d)^p]_{\mathrm{loc}}\virt$ (which are, by Theorem \ref{GW-GSW},
$N_{g,d}$'s up to a sign) are related to the FJRW invariants $\Theta_{g,k}$'s defined by $[M_{g,2^k}^{1/5,5p}]\virt_{\mathrm{loc}}$. 
We understand that  the task is to construct theories that quantitatively link all in
  \beq\label{phases}
 \text{GW of } \ Q \quad  \overset{\textup{Theorem \ref{GW-GSW}}}{\Longleftrightarrow}
 \text{GSW of }\ (K_\Pf,\bw) \quad \longleftarrow \scalebox{1.3}{?} \longrightarrow  \quad
 \text{FJRW of }\ ([\CC^5/\ZZ_5],W).
 \eeq
 
To pursue this goal, we immediately face a specific problem: ``the change of phases' sign". 
The topological type of fields in $K_\Pf$ is labelled by a pair $(g,d)$, where $g$ is the genus  of the curve $C$ and 
$d=\deg L=\deg f\sta \sO(1)$ is always a non-negative integer; the topological type of fields valued in
 $([\CC^5/\ZZ_5],W_5)$ is labelled by a pair $(g,k)$, the genus $g$ of the curve and the number $k$ of $2/5$ marking. 
When $g$ is fixed and when $k$ is general (large) enough, one can show the degree of the line bundle (over the coarse curve) can be arbitrarily negative,
namely, in the phase  $([\CC^5/\ZZ_5],W_5)$, fields are generally of negative degree. In GLSM this corresponds to the fact that $N_{g,d}$ are invariants near
large radius limit point $r\gg 0$ and the LG phase $([\CC^5/\ZZ_5],W_5)$ occurs near the orbifold point $r\ll 0$ (c.f. \cite[Section 5.1]{GLSM}).
 
How can a field of positive degree be transformed to a field of negative degree? 
In which space could this unusual transform happen? How does such transformation -- if it exists -- change the virtual cycles and counting?
We will address these questions in the following sections.


 \section{Master space}
\subsection{Mathematics}  If one builds a large moduli space containing $\barM_{g}(\Pf,d)^p$   
and $\barM_g^{1/5,5p}$ as its \textit{disjoint} closed subspace, then intersection theory over the large moduli 
would give us information relating $\barM_{g}(\Pf,d)^p$  to $\barM_g^{1/5,5p} $. Recall that the two target spaces
$K_\Pf$ and $[\CC^5/\ZZ_5] $ are both open 
subsets of the 5-dimensional stack $[\CC^6/\CC\sta]$, where the two overlap on a large open set
$K_\Pf-\Pf=[(\CC^5-  \{\vec 0\}) /\ZZ_5]$.  If one embeds these 
two open substacks as disjoint closed substacks of a higher dimensional stack $W$, we may
consider the space of maps from curves to $W$ as just stated.
This higher dimensional stack has a natural construction in various places in  
``Variation of GIT" (VGIT) before, called the {\em master space} after M. Thaddeus. Here is a brief introduction. 

Consider the following $\mathbb C^*$-action on $\mathbb C^5\times\mathbb C\times \mathbb P^1$: for $t\in \mathbb C^*$, 
$$
(x_1, \ldots, x_5, p,[u_1, u_2])^t\colon =(tx_1,\ldots, tx_5, t^{-5}p, [tu_1, u_2]).
$$
There is a GIT quotient 
$$
\bar{W} \colon = (\mathbb C^5\times \mathbb  C\times \mathbb P^1-\sS)/\mathbb C^*\hbox{ where }  \sS\colon=\{(x_i=0=u_1)\cup (p=0=u_2)\}
$$
which is a 6-dimensional simplicial toric variety. So $\bar W$ has at most orbifold singularities. Indeed,
$\bar{W}$ is smooth outside the unique orbifold point given by $x_i=u_2=0$.  
The stacky quotient
$$
W =[(\mathbb C^5\times \mathbb  C\times \mathbb P^1-\sS)/\mathbb C^*]
$$ 
is a 6-dimensional smooth toric Deligne-Mumford (DM) stack with coarse moduli space $\bar W$. 

Consider a $\mathbb C^*$-action on $W$, and call this action $T$-action to avoid confusion. For $t\in T=\mathbb C^*$,
$$
(x_1,\ldots, x_5, p, [u_1, u_2])^t=(x_1,\ldots, x_5, p, [tu_1, u_2]).
$$
The $T$-fixed locus is a disjoint union of three connected components:
$$
W^T=K_{\Pf}{\times} \{0\}\sqcup \vec 0\times 
 \left( (\mathbb P^1-\{0,\infty\})/{\mathbb C}^*\right)   \sqcup 
\mathbb [\mathbb C^5/\mathbb Z_5]\times \{\infty\}
$$
where $0=[1, 0]$ and $\infty=[1, 0]$ in $\mathbb P^1$, and the middle term $ \left( (\mathbb P^1-\{0,\infty\})/{\mathbb C}^*\right)  $ is nothing but one single point.
 \begin{figure}
 \includegraphics[scale=0.4]{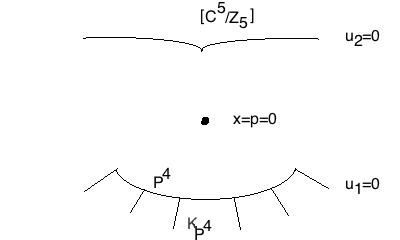}
 \caption{picture for the master space $W$, a six dimensional manifold.}
 \end{figure}
 The shape of $W$ and its $T$ fixed loci is shown in Figure 2, where  
$\sO_\Pf(-5)=K_\Pf$ and $[\CC^5/\ZZ_5]$ are disjoint divisors defined by
$u_1=0$ and $u_2=0$ respectively.
The single point defined by $x_i=p=0$ is responsible for the conifold point of complexified K\"{a}hler moduli space
of the quintic.

\section{Mixed Spin Fields: Quantization of the Master Space }

\subsection{Mixed Spin P-fields}

 Following the recipe from previous sections, now we consider a field theory valued in $W$, namely 
the space of maps from curves to the master space $W$.

Such an objet is called  avmixed spin $P$-field (MSP for short). It consists of
$$
\xi=(\Sigma^\sC, \sC, \sL , \sN, \varphi_1,\ldots, \varphi_5, \rho, \nu=[\nu_1,\nu_2])
$$
where 
\begin{enumerate}
\item $(\Sigma^{\sC}, \sC)$ is a pointed twisted curve, 
\item $\sL$ and $\sN$ are invertible sheaves on $\sC$ ($\sL$ is as before but $\sN$ is new due to the extra factor 
$\mathbb P^1$ in the master space technique),
\item $\varphi_i\in H^0(\sL)$ and $\rho\in H^0(\sL^{\vee 5}\otimes \omega_{\sC}^{\log})$ (as before),
\item $\nu=[\nu_1,\nu_2]$ is a new field, where $\nu_1\in H^0(\sL\otimes \sN)$ and $\nu_2\in H^0(\sN)$.
\end{enumerate}
They satisfy the following conditions: 
\begin{enumerate}
\item (narrow condition) $\varphi_i|_{\Sigma^{\sC}}=0$,
\item (combined GIT-like stability conditions) 
\begin{enumerate}
\item $(\varphi_1, \ldots, \varphi_5, \nu_1)$ is nowhere vanishing (coming from excluding 
$\{(x_i=0=u_1)\}$),
\item $(\rho, \nu_2)$ is nowhere vanishing (coming from excluding $\{(p=0=u_2)\}$),
\item $(\nu_1, \nu_2)$ is nowhere vanishing (coming from $[u_1, u_2]\in \mathbb P^1$). 
\end{enumerate}

\end{enumerate}
We say $\xi$ is stable if Aut($\xi$)  is finite. 
For simplicity, we will use $\varphi$ to represent $(\varphi_1, \ldots, \varphi_5)$.

In order to understand why the moduli space of MSP fields geometrically contains 
the moduli space $\barM_g(\Pf,d)^p$ of stable maps with P-fields and
the moduli space $\barM_g^{1/5,5p}$ of 5-spin curves with five P-fields, we examine the moduli space of MSP fields in details. 

Let $\xi$ be a MSP field. 
\begin{enumerate}
\item When $\nu_1=0$, since $(\varphi_1, \ldots, \varphi_5,\nu_1)$ is nowhere zero, we must have $(\varphi_1, \ldots, \varphi_5)$ is nowhere zero. Since $(\nu_1, \nu_2)$ is nowhere zero, $\nu_2$ must be nowhere zero. Since $\nu_2$ is a section of $\sN$, $\sN\cong \sO_\sC$. There is no restriction on $\rho$. Thus $\xi\in \barM_g(\Pf,d)^p$ and we get GW theory of the quintic threefold $Q$.
\item When $\nu_2=0$, since $(\rho, \nu_2)$ is nowhere zero, $\rho$ must be nowhere vanishing. 
Since $\rho$ is a section of $\sL^{\vee 5}\otimes \omega_\sC^{\log}$, we must have $\sL^5\cong \omega_\sC^{\log}$. 
Also $\nu_1$ must be nowhere zero. Thus $\sL\otimes \sN\cong\sO_\sC$, i.e., $\sN\cong \sL^\vee$.  
$\varphi_1,\ldots, \varphi_5$ can be arbitrary. Thus $\xi\in \barM_{g,(\gamma_j)}^{1/5,5p}$ and we get the FJRW theory of $([\CC^5/\ZZ_5],W_5)$.

\item When $\rho=0$ and $\varphi_i=0$ for $1\le i\le 5$, $\nu_1, \nu_2$ must be nowhere zero. 
Thus $\sN\cong \sO_\sC$ and $\sL\cong\sO_\sC$. Hence we get moduli of stable curves
$\barM_{g,n}$. 
\end{enumerate}

\begin{theo}[H.L. Chang - J. Li - W.P. Li - C.C. Liu \cite{CLLL}] The moduli stack $W_{g, \gamma, {\bf d}}$ of stable MSP fields of genus $g$, monodromy $\gamma=(\gamma_1, \ldots,\gamma_\ell)$ of $\sL$ along $\Sigma^{\sC}$ and degree ${\bf d}=(d_0, d_\infty)$ of $\sL\otimes \sN$ and $\sN$ respectively is a separated DM stack of locally finite type.
\end{theo}

The moduli stack $W_{g, \gamma, {\bf d}}$ admits a natural $\mathbb C^*$-action also called $T$-action: for $t\in \mathbb C^*$, 
\begin{eqnarray*}
(\Sigma^{\sC}, \sC, \sL, \sN, \varphi, \rho, \nu_1, \nu_2)^t\colon =(\Sigma^{\sC}, \sC, \sL, \sN, \varphi, \rho, t\nu_1, \nu_2).
\end{eqnarray*}
It is not  proper since $\varphi$ and $\rho$ are sections of invertible sheaves. Thus we cannot do integrations on this stack. However, there exists a cosection of its obstruction sheaf. Using the arguments similar to the GW case and LG case, we have the following theorem.

\begin{theo}[\cite{CLLL}] The moduli stack $W_{g, \gamma, {\bf d}}$ has a $T$-equivariant perfect obstruction theory, a $T$-equivariant cosection $\sigma$ of its obstruction sheaf, and thus carries a $T$-equivariant cosection localized virtual cycle
$$[W_{g, \gamma, {\bf d}}]\virt_{\mathrm{loc}}\in A_*^T W_{g,\gamma,{\bf d}}^-$$
where $W_{g,\gamma,{\bf d}}^-$ is the degeneracy locus of $\sigma$, i.e.,
$$W_{g,\gamma,{\bf d}}^-\colon=(\sigma=0)=\{\xi\in W_{g,\gamma,{\bf d}}\, |\, \sC=(\varphi=0)\cup (\varphi_1^5+\ldots+\varphi_5^5=0=\rho)\}.$$
\end{theo}

 The cycle  enumerates ``maps to the master space $W$". Figure 3 is an example where the domain curve is represented as a union of one dimension lines (which is the standard notation in algebraic 
geometry).

 \begin{figure}
 \begin{center} 
 \includegraphics[scale=0.5]{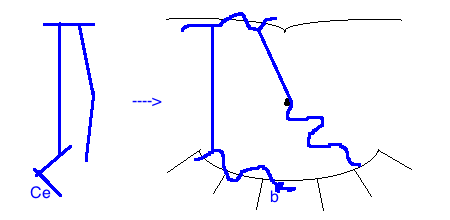}
 \end{center}
 \caption{graphs for a field valued in master space $W$.}
  \end{figure}

  In Figure 3, the component $C_e$ is mapped to single point $b$, and is what we call
a ``ghost'' (over which $\rho$ can be nonvanishing) in Figure 1. Note that considering
 $\xi$ as a map is just for easiness of understanding: indeed the map cannot be realized
 due to the presence of $\omega_\sC$ in the definition of ($\rho$ in) $\xi$.

\subsection{Properness: Capture Ghost at Infinity}

In order to integrate, we need properness of $W_{g,\gamma,{\bf d}}^-$. In fact, we have
\begin{theo}[\cite{CLLL}] \label{thm:proper} The degeneracy locus $W_{g,\gamma,{\bf d}}^-$ is a proper $T$-DM stack of finite type.
\end{theo}

The proof of the  properness reveals an important phenomenon transforming fields of different phases in the MSP moduli. 
Under the transformation, the spin structure of line bundles arises naturally in the LG-phase as a limit of a family of P-fields in CY-phase. 
We call this phenomenon the  ``Landau-Ginzburg transition".  As mentioned before, the contribution from ghosts is one of the difficulties to approach 
postive genus Gromov-Witten invariants. The LG transition phenomenon
enables the FJRW theory to capture the ghost contribution in the GW theory, inside the MSP moduli space.

\subsubsection{LG-Transition: An Example} \label{LG-Transition1}
For any positive integer $d$, we construct a simple example where $g=1$, $\gamma=\emptyset$, and $\bd=(d,0)$
to illustrate the phenomenon of LG-transition and explain why FJRW theory  comes into the picture naturally when we consider GW theory with a P-field. 
The argument below is also a part of the procedure to prove Theorem \ref{thm:proper} (properness of the degeneracy locus).

\noindent
{\it 1. A  point in the degeneracy locus $W_{1,\emptyset, (d,0)}^-$.}

 We give an MSP-field which looks like the picture in the left of Figure 4.
Given a positive integer $d$, define an MSP-field 
\begin{eqnarray}\label{xi}
\xi=(\sC,\sL,\sN,\varphi,\rho,\nu= [\nu_1,\nu_2])
\end{eqnarray}
over a point as follows. The curve $\sC$ is a union of a smooth elliptic curve $C_1$ and a smooth rational curve $C_0$, intersecting at a node $p$.  
Under the isomorphism $C_0\cong \PP^1$ we have
$$
\sL|_{C_0}\cong \cO_{\PP^1}(d),\quad \cN_{C_0}\cong \cO_{\PP^1}, \quad
\varphi|_{C_0}=(x^d, -x^d,0,0,0),\quad\rho|_{C_0}=0,\quad 
\nu_1|_{C_0}=y^d,\quad \nu_2|_{C_0}=1
$$
where $[x,y]$ are homogeneous coordinates on $C_0=\PP^1$, and $p=[0,1]$. In particular
$\varphi(p)=0$ and $\nu_1(p)=\nu_2(p)=1$. On $C_1$, we have 
$$
\sL|_{C_1}\cong \cO_{C_1},\quad \sN|_{C_1}\cong \cO_{C_1}, \quad
\varphi|_{C_1}=(0,\ldots,0),\quad \nu_1|_{C_1}=1,\quad \nu_2|_{C_1}=1.
$$
In particular, $\cN\cong \cO_{\sC}$.
Finally, we extend $\rho|_{C_0}=0$ to a non-zero section $\rho\in H^0(\sL^{\vee 5}\otimes \omega^{\log}_C)$ as follows: $\rho|_{C_1}$ is a non-zero section 
of $H^0\big((\sL^{\vee 5}\otimes \omega_C)|_{C_1}\big)=H^0(\omega_{C_1}(p))$ vanishing at $p$ only.  The choice of $\rho|_{C_1}$ is unique up to
multiplication by a nonzero constant. Then $\xi$ represents a point in 
the degeneracy locus $W_{g=1,\gamma=\emptyset,\bd=(d,0)}^- \subset W_{g=1,\gamma=\emptyset, \bd=(d,0)}$. 

\medskip

\noindent
{\it 2. A morphism from $\CC^*$ to $W_{1,\emptyset, (d,0)}^-$}.

 We describe a one-parameter deformation of the MSP field $\xi$, depicted by Figure 4.
Let $S_* = \sC\times \CC^*$ and let $\pi_1: S_*\to \sC$ be the projection to the first factor.
We consider a family of MSP-fields over $\mathbb C^*$
$$
\xi_*=(S\lsta ,  \cL\lsta=\pi_1^*\sL, \cN\lsta=\pi_1^*\sN =\cO_{S\lsta}, 
\varphi\lsta=\pi_1^* \varphi, \rho\lsta=t^{-1}\pi_1^*\rho, [\nu_{1*}=\pi_1^*\nu_1, \nu_{2*}=\pi_1^*\nu_2=1])
$$
where $t$ is the parameter of $\CC^*=\mathrm{Spec}[t,t^{-1}]$.  This family over $\CC^*$ defines a morphism
\begin{equation}\label{eqn:Cstar-to-W}
\phi: \CC^* \longrightarrow W_{1,\emptyset, (d,0)}^-. 
\end{equation}

By abuse of notation, let $\pi_1$ also denote the projection
from $S_{*i}:= C_i\times \CC^*$ to the first factor, where $i=0,1$.
The restriction of the family $\xi_*$ to 
$S_{0*}$ is a constant family over $\CC^*$:
$$
\xi_{0*}= (\Si_{0*}= p\times \CC^*, S_{0*}= C_0\times \CC^*,
\pi_1^*\cO_{\PP^1}(d), \cO_{S_{0*}}, (x^d, -x^d, 0,0,0), 
0, [y^d,1])
$$
which defines a {\em constant} map 
$\phi_0:  \CC^*\longrightarrow W_{g=0, \gamma=(1), \bd=(d,0)}^-$.
The restriction of the family $\xi_*$ to  $S_{1*}$ is
$$
\xi_{1*}= (\Si_{1*}= p\times \CC^*, S_{1*}= C_1\times \CC^*,
\cO_{S_{1*}}, \cO_{S_{1*}}, (0,\ldots,0), 
t^{-1}\pi_1^*(\rho|_{C_1}), [1,1])
$$
which defines a morphism $\phi_1: \CC^*\longrightarrow 
W_{g=1,\gamma=(1),\bd=(0,0)}$.

\medskip

\noindent
{\it 3. The limits $t\to 0$ and $t\to \infty$}\\
We will see that the morphism \eqref{eqn:Cstar-to-W} extends to a morphism
$$
\bar{\phi}: \PP[1,5] \longrightarrow W_{1,\emptyset,(d,0)}^-
$$
where the embedding $\CC^*\hookrightarrow \PP[1,5]$ is given by
$t\mapsto [1,t^{-1}]$. The image $\bar{\phi}([0,1])$ (resp. $\bar{\phi}([1,0])$) is the limit in $W_{1,\emptyset,(d,0)}^-$  when $t\to 0$ (resp. $t\to \infty$). It is easy to see that
$$
\bar{\phi}([1,0]) = (\sC, \sL, \sN, \varphi, 0, \nu=[\nu_1 ,\nu_2])
$$
where $\sC, \sL, \sN, \varphi, \nu_1,\nu_2$ are defined as in Step {\it 1}.

 \begin{figure}[H]
 	\includegraphics[scale=0.5]{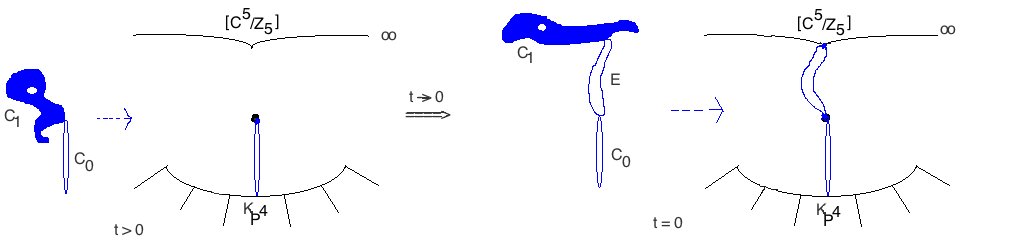}
 	\caption{via scaling P-field $\rho\mapsto \rho/t $ and making $t\to 0$, the ghost is capture  at   infinity  and become FJRW instantons.}
 \end{figure}

 The limit becomes a new MSP field which is the picture in the right of Figure 4. 
 Below we provide detailed construction of the family, which may be technical.\\

The extension of the constant map $\phi_0: \CC^*\longrightarrow W_{0,\gamma,\bd=(d,0)}^-$ is the constant map $\bar{\phi}_0: \PP[1,5]\to W_{0,\gamma,\bd=(d,0)}^-$. To find the limit $\bar{\phi}([0,1])$, it suffices to find $\bar{\phi}_1([0,1])$ where $\bar{\phi}_1$ is the extension of $\phi_1$ to  $\PP[1,5]$.

Let $S_1=C_1\times \mathbb C^*$ and $\rho_1=\rho|_{C_1}$. If we naively take $\cL=\cO_{S_1}$, $\cN=\cO_{S_1}$, $\varphi_{i}=0$, , $\nu_{1}=1$ and $\nu_{2}=1$,  the section $ t^{-1}\pi_2^*\rho_1$
cannot be extended to a regular section of $\sL^{\vee5}\otimes \omega^{\log}_{S_1/\mathbb C}$. 
 Here by abuse of notations, $\pi_i$ is the projection from $S_1$ to the $i$-th factor. One way to solve this problem is to use the equivalence of $\xi_{1*}$ with the following $\xi_{1*}^\prime$:
  $$(S_{1*},  \Sigma_{1*}, \cL_{ *}=\cO_{S_{1*}}, \cN\lsta=\cO_{S_{1*}}, \varphi_{1*}=0, \ldots,  \varphi_{5*}=0, \rho_{S_1*}=\pi_1^*\rho_1, \nu_{1*}=1, \nu_{2*}=t^{\frac{1}{5}}).$$
  Then we can have the extension 
    $$(S_1,  \Sigma_{1}, \cL=\cO_{S_1}, \cN=\cO_{S_1}, \varphi_{1}=0, \ldots,  \varphi_{5}=0, \rho_{S_1}=\pi_1^*\rho_1, \nu_{1}=1, \nu_{2}=t^{\frac{1}{5}})$$
    where $\Sigma_1=p\times \mathbb C$.  
    
The term $t^{\frac{1}{5}}$ may look troublesome. Let's just treat this as indicating  that the zero locus of $\nu_2$ is $\frac{1}{5} C_1\times 0$. The issue of fractional divisor will be resolved once we work in the world of twisted curves. 

Let's assume that we can work with fractional divisors. The zero divisor of $\rho_{S_1}$ is $\Sigma_1$ and that of $\nu_{2}$ is $\frac{1}{5}C_1\times 0$. Since these two divisors intersect,  by MSP requirement that $\rho_{S_1}$ and $\nu_2$ cannot  be zero simultaneously at any point , we don't get an MSP extension. Thus we need to blow up the intersection of these two divisors to separate them. 

Let $\tau\colon S_1^\prime\to S_1$ be the blow up of $S_1$ at $P=\Sigma_{1}\cap (C_1\times 0)=p\times 0$,  $E$ be the exceptional divisor, $\Sigma_1^\prime\subset S^\prime $ be the strict transform of $\Sigma_{1}$, and $C_1^\prime$ be the strict transform of $C_{1}\times 0$. The zero divisor of $\tau^*\nu_{2}$ is $\frac{1}{5}C_1^\prime+\frac{1}{5}E$. 
The zero divisor of $\tau^*\rho_{S_1}$ is $\Sigma_1^\prime+ E$. Now we need to modify $\cL$ and $\cN$ by replacing $\cL$ by $\cL^\prime=\cO_{S_1^\prime}(E/5)$ and $\cN$ by $\cN^\prime=\cO_{S_1^\prime}(-E/5)$. Here we pretend that $\cL^\prime$ and $\cN^\prime$ exist. Indeed, they do not exist in the ordinary sense, but their existence  again will be resolved once we work with orbifolds.  Let $\rho^\prime$ be the section in $H^0(\cL^{\prime\vee 5}\otimes \omega^{\log}_{S_1^\prime/\mathbb C})=H^0(\omega^{\log}_{S_1^\prime/\mathbb C}(-E))$ whose image is $\tau^*\rho_{S_1}$ under the natural map $H^0(\omega^{\log}_{S_1^\prime/\mathbb C}(-E))\to H^0(\omega^{\log}_{S_1^\prime/\mathbb C})$. The zero divisor of $\rho^\prime$ is $\Sigma_1^\prime$.  Let $\nu_{2}^\prime$ be the section of $\cO_{S_1^\prime}(-E/5)$ whose image under the natural map  $\cO_{S_1^\prime}(-E/5)\to \cO_{S_1^\prime}$ is $\tau^*\nu_{2}$. Then the zero divisor of $\nu_2^\prime$ is $\frac{1}{5}C_1^\prime$. Since $C_1^\prime$ and $\Sigma_1^\prime$ don't intersect, we can get an MSP extension by taking $\cL^\prime=\sO_{S_1^\prime}(E/5)$, $\cN^\prime\cong \cL^{\prime\vee}$, $\nu_1^\prime$  a nonzero constant section  $\sL^\prime\otimes \sN^\prime$, $\varphi_1^\prime=\ldots=\varphi_5^\prime=0$, $\rho^\prime\in H^0(\cL^{\prime\vee 5}\otimes \omega^{\log}_{S_1^\prime/\mathbb C})$ with its zero divisor being $\Sigma_1^\prime $ which is the marking of $S_1^\prime /\mathbb C$, and $\nu_2^\prime\in H^0(\cN^\prime)$ whose zero divisor is $\frac{1}{5}C_1^\prime$.  As we mentioned earlier, to make this construction rigorous, we have to do base change and introduce stacky structures at nodes of $C_1^\prime\cup E$ (see \cite{AGV, CLLL}).

Now we see that  the central fiber of $S_1^\prime $ over $\mathbb C$ at $0\in \mathbb C$ is set-theoretically  $ C_1\cup E$, a union of the elliptic curve $C_1$ with the smooth rational curve $E$. The section $\nu_2^\prime$ vanishes on $\sC_{\infty}=C_1^\prime\cong C_1$ and $\rho^\prime$ is nowhere vanishing on $\sC_\infty$. Hence $\sC_\infty$ is a $5$-spin twisted curve, and $E$ is a rational smooth  twisted curve with a marking where $\rho^\prime$ vanishes. 

Then we can glue the MSP field on $S_1^\prime$ with the MSP field on $C_0\times\mathbb C$ by identifying  the marking $\Sigma_1^\prime$ with the marking $p\times \mathbb C\in C_0\times \mathbb C$ after  possibly a base change. Thus the central fiber of the extension is a union of a smooth rational curve $C_0$,  an elliptic curve $C_1$ which is a $5$-spin curve, and a rational twisted curve $E$    intersecting with $C_1$ at the stacky point and with $C_0$ at another point where the nonzero $\rho$-section vanishes.  

 We can also deform the MSP field \eqref{xi} to a MSP field in GW sector as follows. 

Consider $S=\sC\times \mathbb C$, Let $\pi_1$ be the projection of $S$ to its first factor. Take,  for $t\in \mathbb C$, we have a family of MSP over $\mathbb C$,
$$
(S, \quad \pi_1^*\sL,\quad\pi_1^*\sN, \quad
\varphi_S,\quad \pi_1^*\rho,\quad
[t\pi_1^*\nu_1,\pi_1^*\nu_2]).
$$
Here $\varphi_S$ is defined as follows. 
$$\varphi_S|_{C_0\times \mathbb C}=(x^d, -x^d, (1-t)y^d, -(1-t)y^d, 0),\quad \varphi_S|_{C_1\times \mathbb C}=(0, 0, 1-t, t-1, 0).$$
When $t=1$, we get the MSP field $\xi$ in \eqref{xi}. When $t=0$, we get an MSP lying in GW sector since $t\pi_1^*\nu_1|_{t=0}=0$.
\black 
\qed

\section{Vanishing and Polynomial Relations}

How to extract information of GW  and/or FJRW invariants from the cycle $[W_{g, \gamma, \bd}]\virt_{\mathrm{loc}}$? 
In this section, we consider a less general case $\gamma=\emptyset$ (i.e. no markings) to illustrate  the key ideas. 
By virtual dimension counting, we have
$$
[W_{g, \bd}]\virt_{\mathrm{loc}}\in H^{\mathbb C^*}_{2(d_0+d_\infty+1-g)}(W_{g, \bf d}^-, \mathbb Q).
$$
When $d_0+d_{\infty}+1-g>0$, letting $u=c_1({\bf 1}|_{wt=1})$, i.e. $u$ is the parameter for $H_{\mathbb C^*}^*(pt)$, we have
\begin{eqnarray*}
[u^{d_0+d_\infty+1-g}\cdot [W_{g, \bd}]\virt_{\mathrm{loc}}]_0=0.
\end{eqnarray*}
Here $[\cdot]_0$ is the degree zero term in the variable $u$. 

Let $\Gamma$ be a graph associated to fixed points of the $T$-action of $W_{g, \bd}$ and $F_{\Gamma}$ be a connected component 
of $W_{g, \bd}^{T}$ of the graph type $\Gamma$. Applying the cosection localized version  of the virtual localization formula of Graber-Pandaripande \cite{GP} proved 
by Chang-Kiem-J.Li in \cite{CKL}, we obtain
\begin{eqnarray}\label{virtual-local}
\sum_{\Gamma}\left[u^{d_0+d_\infty+1-g}\frac{[F_\Gamma]^{vir}_{loc}}{e(N_{F_\Gamma})}\right]_0=0.
\end{eqnarray}
To deal with $[F_\Gamma]\virt_{\mathrm{loc}}$, we need a decomposition result to be explained below.

Let $\xi=(\sC,  \sL, \sN, \varphi, \rho, \nu_1,\nu_2) \in (W_{g, \bd})^T$ be an MSP field fixed by the $T$-action. We set
\begin{enumerate}
\item $\sC_0$ to be the part of $\sC$ where $\nu_1=0$;
\item $\sC_1$ to be the part of $\sC$ where $\varphi=0=\rho$ and hence $\nu_1=1 =\nu_2$, i.e., $\nu_1$ and $\nu_2$ are nowhere zero;
\item $\sC_\infty$ to be the part of $\sC$ where $\nu_2=0$.
\end{enumerate}
Thus 
\begin{itemize}
\item $\xi|_{ \hbox{ (connected component of $\sC_0$)} }$ is in $\barM_{g^\prime, n^\prime}(\Pf, d^\prime)^p$ which gives GW invariants of $Q$. 
Here marked points appear coming from some nodes on $\sC_0$. 
\item $\xi|_{\hbox{(connected component of $\sC_1$})}$ is in $\barM_{g^\prime, n^\prime}$ which gives Hodge integrals.
\item  $\xi|_{\hbox{(connected component of $\sC_\infty$})}$ is in $\barM^{\frac{1}{5}, 5p}_{g^\prime, \gamma^\prime}$ which gives FJRW invariants of $(\CC^5/\ZZ_5, W_5)$ where $\gamma^\prime$ appears because of some stacky nodes on $\sC_\infty$. 
\end{itemize}

We have the following decomposition result:
\begin{eqnarray*}
[F_\Gamma]^{vir}_{loc}=c\prod [\hbox{moduli of $\xi|_{\sC_0}$]}\virt_{\mathrm{loc}}
\cdot[\hbox{moduli of $\xi|_{\sC_1}$}]\virt_{\mathrm{loc}}\cdot [ \hbox{moduli of $\xi|_{\sC_\infty}$}]\virt_{\mathrm{loc}}
\end{eqnarray*}
where $c$ is a constant.  The first factor gives GW invariants of stable maps to $\PP^4$ with P-fields, i.e. $N_{g',d'}$. 
The second  factor gives Hodge integrals on $\overline M_{g',n'}$. The third factor gives FJRW invariants of insertions $-\frac{2}{5}$ (after using a vanishing).  After $e(N_{F_\Gamma})$'s are calculated, using the polynomial relations (\ref{virtual-local}), we obtain the following results about GW invariants of the quintic.

\begin{theo}[\cite{CLLL2}]\label{thm-induction}
Letting $d_\infty=0$, the relations \eqref{virtual-local} provide an effective algorithm to evaluate the GW invariants $N_{g,d}$
provided the following are known
\begin{enumerate}
\item $N_{g',d'}$ for $(g',d')$ such that $g'< g$, and $d'\le d$;
\item $N_{g,d'}$ for $d'< g$;
\item $\Theta_{g',k}$ for $g'\leq g-1$ and $k\leq 2g-4$;
\item $\Theta_{g,k}$ for $k\leq 2g-2$.
\end{enumerate}
\end{theo}
 
Recall that $\Theta_{g,k}$ is the genus $g$ FJRW invariants of insertions $-\frac{2}{5}$ and $\Theta_{g,k}$ may be non-zero only when $k+2-2g\equiv 0 (5)$.
We can see that when $g=2$ only  $\Theta_{2,2}$ is needed,  and when $g=3$ only $\Theta_{3,4}$ is needed.  

\begin{rema}  As we know,  on using mathematical induction, upon more numerical datum  the induction is,  the less effective the computation will be.  
We can see from Theorem \ref{thm-indcution} that MSP induction for GW invariants is carried out on two numbers, genus and the degree only. 
Thus this provides a rather effective way to facilitate the induction procedure.
\end{rema}

We can also use the vanishings \eqref{virtual-local} for $d=(0,d_\infty)$,  to determine quintic's FJRW invariants up to finite many initial data.

\begin{theo}[\cite{CLLL2}]\label{nomarking}
 For a fixed positive genus $g$, the finite set $\{\Theta_{g,k}\}_{k<7g-2}$ determine all genus $g$ FJRW invariants  $\{\Theta_{g,k}\}_{k=0}^\infty$.
\end{theo} 
These relations are effective in calculating FJRW invariants. For example, for the case of genus $2$, $\{\Theta_{2,k}\}_k$ can be inductively derived from only two unknowns $\Theta_{2,2}$ and $\Theta_{2,7}$.

We end this section by some speculations.

 Let us look at Theorem \ref{thm-induction} from a different aspect. Inductively we may suppose all GW/FJRW invariants for genus less than $g$ are known. Then for genus $g$, Theorem \ref{thm-induction} reduces the problem of determining  the infinitely  many  GW invariants  $\{N_{g,d}\}_{d=1}^\infty$ to two finite sets of initial datum
 $$ \{N_{g,1},\cdots,N_{g,g-1} \}  \and  \{\Theta_{g,k}\}_{k\leq 2g-2}.$$


We formulate the following speculation:

\medskip

\textit{By suitable choice of positive $d_0$ and $d_\infty$, the relations (\ref{virtual-local}) provide an effective algorithm to determine the first set of initial data $ \{N_{g,1},\cdots,N_{g,g-1} \} $.}

\medskip

 If this is true,  then one is left to determine the second set of initial data $  \{\Theta_{g,k}\}_{k\leq 2g-2}$.  We propose another conjecture about fully determining all FJRW invariants for the quintic,

\begin{conj}\label{conj-2}
The  equations  (\ref{virtual-local})  using $d_0=0$ 
and nonempty $\gamma$'s (i.e. with markings)  give relations  that effectively evaluate all $\Theta_{g,k}$.
\end{conj}

\section{Comparison with Physical Theories}

\subsection{Comparison with Witten's GLSM}

In \cite{GLSM}, Witten introduced a family of theories using path integrals, 
called the Gauged Linear Sigma Model (GLSM), linking a non-linear sigma models on a Calabi-Yau hypersurface to a
sigma model targeting  in a Landau-Ginzburg space.
The GLSM is parameterized by a complexified K\"{a}hler parameter 
$$t:=r-i\theta$$
where  $r$ is  ``Fayet-Iliopoulos parameter"  and $\theta$ is called the theta angle\footnote{we use notations in \cite[Sect 15.2.2]{MS};}. 
Witten \cite[Sect 3.1]{GLSM} argued that the GLSM specializes to GW path integral when $r\to +\infty$,
and specializes to LG model path intergral when $r\to -\infty$. This is known as the
Calabi-Yau/Landau-Ginzburg correspondence.

The Mixed-Spin-P fields (MSP fields) introduced in \cite{CLLL}  is  a field theory designed to capture ``phase space
transition" in one cage\footnote{in the MSP cage (proper) integral does not diverge because the cage is proper and separated by \cite{CLLL}}. 
 An MSP field can be viewed as an interpolation between
fields valued in $K_\Pf$ and fields valued in $[\CC^5/\ZZ_5]$, and the interpolation is governed by the  
``$\nu$ field". Over the part of worldsheet (curve) where $\nu=0$, the MSP field is a pure field taking values in $K_{\Pf}$, and, over
$\nu=\infty$, the MSP field is a pure field taking values in $[\CC^5/\ZZ_5]$.
In a nutshell, by promoting the phase parameter $r-i\theta$ into a field $\nu$ on worldsheet (curve), 
we transform Witten's family of theories parametrized
by $r-i\theta$ into a single new field theory.  Also, an advantage of MSP moduli is it works for higher loop (in physics terms)
or for higher genus (in math terms), while GLSM in physics literature does not\footnote{GLSM only treat genus zero worldsheet}. 

 We recall the question raised by Witten \cite[Page 28]{GLSM}:
``\textit{Are Calabi-Yau and Landau-Ginzburg separated by a true phase transition, at or near $r = 0$?
There is no reason that the answer to this question has to be �universal,� that is, independent of the path one follows in interpolating from Calabi-Yau to Landau-Ginzburg in a multiparameter space of not necessarily conformally invariant theories. Along a suitable path, there may well be a sharply defined phase transition, while along another path there might not be one. This seems quite plausible.}"
Though our $\nu$ field does not have a definite value of which we can vary ``the theories",  $\nu$ allows
us to introduce a $T=\CC\sta$ action, and by localizing to $T$ fixed locus, we obtain (possibly as referred to by Witten) a
precise phase transition between the quintic Calabi-Yau threefold  and the LG of the Fermat quintic.
Furthermore, such phase transitions are multi-fold: for each $\bd=(d_0,d_\infty)$ that provides a vanishing, we get 
an interpolation. And when we vary  $\bd$, we obtain a class of such interpolations.  Thus in physics terms, we may say that 
each $(d_0,d_\infty)\in \ZZ_{\geq 0}\times \ZZ$, being a mode describing how worldsheet is wrapped to K\"{a}hler moduli spaces, 
provides a path linking CY point to LG point in the phase spaces.

\subsection{Compare with B-model} 
For the quintic Calabi-Yau threefold,  the modularity of the generating function $F^A_g$ is suggested by { physicists, 
but is a mystery in mathematics}. In physics literature as \cite{BCOV,HKQ}, the modularity of $F^B_g$ and the mirror symmetry\footnote{another mystery in mathematics} $F^B_g=F^A_g$  When $g>0$, the holomorphic anomaly equation determines $F_g^B(q)$ up to $3g-2$ unknowns. 
The degree zero Gromov-Witten invariant $N_{g,d=0}$ is known, so we are left with $3g-3$ {\blue unknowns.}
The boundary conditions at the orbifold point (which corresponds to
Landau-Ginzburg theory of the Fermat quintic polynomial in five variables) impose $\lceil\frac{3}{5}(g-1)\rceil$ contraints
on the $3g-3$ unknowns, whereas the ``gap condition'' at the conifold point imposes $2g-2$ constraints on the $3g-3$ unknowns. In summary, the holomorphic anomaly equation and the boundary conditions  determine $F_g^{A(Q)}$ up to $\lfloor \frac{2}{5}(g-1)\rfloor$ unknowns.
Coincidently, granting Conjecture A, the number of initial data needed to determine $F_g$ are the FJRW invariants 
$\varTheta_{g,k\leq 2g-2}$, subject to $2g-2\equiv k(5)$. Thus $\lfloor\frac{2(g-1)}{5}\rfloor+1$  many FJRW invariants are
needed to determine $N_{g,d}$ via MSP moduli, provided all lower genus invariants $\{F_{g'}: g'<g\}$ are known. 
We hope there is more geometric explanation then viewing it just as a coincidence.

\end{document}